\documentclass[leqno,12pt]{amsart}
\usepackage{amssymb,amsmath,amscd}

\newtheorem{thm}{Theorem}[section]

\newtheorem{corollary}[thm]{Corollary}
\newtheorem{prop}[thm]{Proposition}
\newtheorem{lemma}[thm]{Lemma}
\newtheorem{fact}[thm]{Fact}

\theoremstyle{definition}
\newtheorem{defn}[thm]{Definition}
\newtheorem{example}[thm]{Example}

\theoremstyle{remark}
\newtheorem{remark}[thm]{Remark}

\newcommand{\bt}{\begin{thm}}
\newcommand{\et}{\end{thm}}
\newcommand{\bp}{\begin{prop}}
\newcommand{\ep}{\end{prop}}
\newcommand{\bd}{\begin{defn}}
\newcommand{\ed}{\end{defn}}
\newcommand{\bl}{\begin{lemma}}
\newcommand{\el}{\end{lemma}}
\newcommand{\bfa}{\begin{fact}}
\newcommand{\efa}{\end{fact}}
\newcommand{\bc}{\begin{corollary}}
\newcommand{\ec}{\end{corollary}}
\newcommand{\bex}{\begin{example}}
\newcommand{\eex}{\end{example}}
\newcommand{\br}{\begin{remark}}
\newcommand{\er}{\end{remark}}
\newcommand{\ben}{\begin{enumerate}}
\newcommand{\een}{\end{enumerate}}

\newcommand{\rrr}{\rightarrow}
\newcommand{\ra}{\rightarrow}

\newcommand{\exact}[3]
{0 \rrr #1 \rrr #2
\rrr #3 \rrr 0}

\newcommand{\PP}{\mathbb{P}}

\newcommand{\coo}{{\mathcal O}}

\newcommand{\cae}{{\mathcal E}}

\begin{document}

\title{
A curve algebraically but not rationally uniformized by radicals}

\author{Gian Pietro Pirola and Enrico Schlesinger}

\address{Dipartimento di Matematica "F. Casorati", Universit\'{a} di Pavia,
via Ferrata 1, 27100 Pavia,  Italia}


\address{Dipartimento di Matematica, Politecnico di Milano, Piazza Leonardo da
Vinci 32, 20133 Milano, Italia}


\thanks{This paper was written in the framework of the research project "Curve e Monodromia"
financed by GNSAGA.
The first author was partially supported by:
1) MIUR PRIN 2003:
  {\em Spazi di moduli e teoria di Lie};
2) Gnsaga; 3) Far 2002 (PV):
  {\em Variet\`{a} algebriche, calcolo algebrico, grafi orientati e topologici.}
The second author was partially supported by
 MIUR PRIN 2002 {\em Geometria e classificazione delle variet\`a proiettive complesse}.
 }

\begin{abstract}
Zariski proved the general complex projective curve of genus $g>6$ is not
rationally uniformized by radicals, that is, admits no map to $\PP^1$ whose Galois group
is solvable.
We give an example of a genus seven complex projective curve $Z$ that is not
rationally uniformized by radicals, but such that there is a finite
covering $Z' \ra Z$ with $Z'$ rationally uniformized by radicals.
The curve providing the example appears in a paper by Debarre and Fahlaoui
where a construction is given to show the Brill Noether loci $W_d(C)$
in the Jacobian of a curve $C$ may contain translates of abelian subvarieties
not arising from maps from $C$ to other curves.

\end{abstract}

\maketitle

\section{Introduction}
Zariski, solving a problem posed by Enriques at the Congress of Mathematicians
held in Zurich in 1897, proves in~\cite{Zariski1}
 that, given an algebraic equation $f(x,y)=0$ of genus $p > 6$ with general moduli,
it is not possible to introduce a parameter $t$, rational function of $x$ and $y$,
in such a way that $x$ and $y$ can be written by radicals as functions of $t$.
Zariski thought of this result as the analogue for algebraic curves of Abel's
Theorem on the non solvability by radicals of a general algebraic
equation of degree $n \geq 5$: in his own words, {\it Questo teorema rappresenta la vera
estensione alle curve algebriche del teorema di Abel}.

We will say - cf.~\cite{fried} and \cite{fg}-  that a complex projective
curve $X$ is {\it rationally uniformized by radicals}
if there is a branched covering $X \ra \PP^1$ whose Galois group is solvable
(by Galois group we mean
the Galois group of the splitting field of the extension of
function fields $\mathbb{C} (X)/\mathbb{C}(\PP^1)$).
Zariski's theorem says that a general smooth complex projective curve of genus
$g > 6$ is not rationally uniformized by radicals. The result is sharp because every curve of
genus $g \leq 6$ arises as a covering of $\PP^1$ of degree $\leq 4$,
hence is rationally uniformized by radicals.

This seminal paper of Zariski has spurred a lot of further
research, especially after the article ~\cite{fg} by Fried and Guralnick. A report on
this work is in \cite{fried}, where among other things a
combinatorial technique is developed to compute the dimension of
the locus $M_g (G)$ of curves with a prescribed Galois group $G$
in the moduli space $M_g$. We
recommend \cite{fried} also for its review of the whole of
Zariski's work on the subject, and for further topics such as
Nielsen classes, branch cycles and their relation to the
endomorphism ring of the Jacobian.

We recall here only  results about Galois
groups in the {\it generic} case, that is, Galois groups of
coverings $f: X \ra \PP^1$ where $X$ is {\it general} in $M_g$.
Any such covering factors as a primitive covering $h:X
\ra \PP^1$ followed by a covering $k: \PP^1 \ra \PP^1$. One is
therefore  led to  study the Galois groups of  primitive coverings
$X \ra \PP^1$ with $X$ general, and  of genus zero coverings.
There have been many developments on the genus zero
problem: see
\cite{gt,fm,fried-genus0,fried-prelude,fried-thompson}.
As for the primitive generic case,
in a series of papers ~\cite{gn,gm,gs}
it is proven that, for any fixed $g \geq 4$,
the Galois group of a  primitive covering $f: X \ra \PP^1$, with $X$
general in $M_g$, is either $A_n$ or $S_n$, with $n >
\frac{g+1}{2}$, a statement that strengthen Zariski's result.

\vspace{0,5cm}

In this paper we are concerned with a fundamental
 question raised by Zariski in the remarks
following its theorem. He writes:
\begin{quote}
"Si potrebbe dunque pensare che si possa invece fornire per {\em ogni}
equazione $ f = 0$ una risoluzione {\em multipla} per radicali $x=x(t)$,
$y=y(t)$, in cui ad ogni punto $(x,y)$ della curva $f=0$ corrispondano
pi\`u valori di $t$. $\ldots$ \`{E} poco probabile che ci\`o accada effettivamente,
ma in ogni modo si ha qui un nuovo problema, che noi non discutiamo in questa Nota
e che potr\`a essere oggetto di una ulteriore ricerca."
\end{quote}
(One may therefore think that for {\em every} equation $f=0$
one can find a {\em multiple}
solution by radicals $x=x(t)$, $y=y(t)$, in which  several values of $t$ correspond
to every given  point $(x,y)$ of the curve $f=0$. $\ldots$
It is unlikely that this could really happen, but in any case we have here
a new problem, which we do not discuss in this Note, and which might be
object of further research).

We define a curve $X$ to be {\em algebraically}
uniformized by radicals if there exists a finite covering $X' \ra X$
where $X'$ is rationally uniformized by radicals.
Zariski is suggesting the general curve $X$ of a fixed, large enough genus should not be
algebraically uniformized by radicals. We do not know
of any progress on this question since Zariski's paper.

In this paper we establish that the two notions of being
rationally,  respectively algebraically, uniformized by radicals are distinct.
We give an example of a curve algebraically, but not rationally,
uniformized by radicals.
The property of not being algebraically uniformized by radicals
is a very strong condition of hyperbolic type:
if $X$ is not algebraically uniformized by radicals and $X^{(4)}$ is the $4$th symmetric product of $X$,
then there are no nonconstant maps $Y \ra X^{(4)}$ for any $Y$ that is algebraically uniformized by
radicals - see \ref{hyp}.
In particular, there are no such maps with $Y$ of genus $g \leq 6$ (or hyperelliptic or trigonal etc.),
and therefore $X^{(4)}$ is hyperbolic, that is, there are no nonconstant maps from an abelian variety to
$X^{(4)}$ - cf.~\cite{ghyp}.

The abelian subvarieties of the symmetric products $X^{(d)}$ have been studied by
Abramovich and Harris in \cite{ah} in connection with the results of Faltings on the
rational points of subvarieties of an abelian variety. In \cite{ah} it is proven that
$X^{(4)}$ is hyperbolic if $X$ is a curve of genus $g \geq 8$ that admits no map
of degree $\leq 4$ onto a curve of genus $\leq 1$ and has no map of degree $2$ onto a
curve of genus $2$.

On the other hand, Debarre and Fahlaoui in \cite{df} give an interesting
example of a curve $Z$ of genus $7$
such that $Z^{(4)}$ contains an elliptic curve, therefore is not hyperbolic,
but $Z$ has no map of degree $\leq 4$ onto a curve of genus $\leq 1$ and no map
onto a curve of genus $g \geq 2$.

In this paper we prove the general curve $Z$ in the family constructed  in \cite{df}
is algebraically but not rationally uniformized by radicals.
$Z$ is algebraically uniformized by radicals because $Z^{(4)}$ contains
an elliptic curve. On the other hand, given the results of \cite{df}, in order
to show $Z$ is not rationally uniformized by radicals it is enough to prove
$Z$ admits no map of degree
$d \geq 5$ onto a curve $W$ of genus $\leq 1$ with {\em primitive and solvable} Galois group.
We prove this counting moduli following the main idea of Zariski's paper. For this
we first need to improve a lemma by Zariski about the possible number of fixed points
of an element $g$ of a primitive solvable subgroup of $S_d$ (Propositions
~\ref{fix1} and~\ref{fix2}).
In section~\ref{three}, following Zariski,
we deduce a lower bound for the branching order at any branch point of a primitive solvable
cover. We thus obtain an upper bound for the dimension of
the family of curves that arise as primitive solvable covers of rational or
elliptic curves.

We think the question of Zariski deserves to be the object of further research,
being a  simple question of  basic importance in the theory of
algebraic curves. We would also like to remark the analogous questions in higher dimension are completely
unexplored. One does not even know whether every variety of dimension $n \geq 2$ admits a dominant map to $\PP^n$
with a solvable Galois group.

We would like to thank M. Freid for many valuable suggestions, and R. Guralnick who
provided us with helpful references and
showed us how to shorten the proof of Proposition~\ref{fix2}.

\vspace{.2in}
\noindent {\bf Notation and terminology.}
We work over the field of complex numbers.
All curves are assumed to be smooth and projective algebraic.
We say that a degree $d$ branched covering
$\pi: X \ra Y$ of two curves is primitive if it does not factor
nontrivially.
This is equivalent to requiring the Galois group $G$ of $\pi$
be a primitive subgroup of the automorphism group $S_d$ of the generic fiber.
We say that $\pi$ is solvable if its Galois group is solvable. We say
a curve $X$ is rationally uniformized by radicals if there exists a solvable
covering $X \ra \PP^1$. We say $X$ is algebraically uniformized by radicals
if there exists a finite covering $X' \ra X$ with $X'$
rationally uniformized by radicals.

\section{Fixed points and primitive solvable groups} \label{two}
Let $\Omega$ be a finite set with $d$ elements, and write
$S_d = \mbox{Aut} (\Omega)$. Let $G \subset S_d$ be a subgroup. Given $x \in \Omega$,
we denote by $G_x$ the stabilizer of $x$ in $G$ and by $Gx$ the $G$-orbit of $x$.
The subgroup $G$ is primitive if it is transitive and $G_x$ is a maximal
subgroup of $G$. Since $G \subset S_d$, the stabilizer $G_x$ contains no proper
normal subgroup of $G$.

The following proposition is well known:
\bp \label{prim}
Let $G \subset S_d$ be a primitive solvable subgroup. Fix $x \in \Omega$. Then
\begin{enumerate}
  \item $G$ contains a unique minimal normal subgroup $A$;
  \item $A$ is an elementary abelian $p$-group for some prime $p$;
  \item $G = A G_x$ and $A \cap G_x = \{1_G\}$;
   \item the action of $A$ on $\Omega$ is regular, that is, for every $x \in \Omega$,
  the map $a \mapsto ax$ is a bijection of $A$ onto $\Omega$.
\end{enumerate}
In particular,  $d = \# A = p^k$ for some $k \geq 1$.
\ep
For a given subset $S$ of $G$, we let
$$ \Omega^S = \{ x \in \Omega: \; gx=x \; \mbox{for every $g\in S$} \}$$
denote the set of fixed points of $S$. We let $G^{\#} = G \smallsetminus \{ 1_G \}$.
\bp[Zariski] \label{fix1}
Let $d=p^k$, and let $G \subset S_d$ be a primitive solvable subgroup.
Then an element $g \in G^{\#}$ has at most $p^{k-1}$ fixed points:
  $
  \# \, \Omega^g \leq p^{k-1}.
  $
\ep
We will need the following refinement:
\bp \label{fix2}
Suppose $d=2^k$ and $d-1$ is a prime (thus a Mersenne number). If
$G \subset S_d$ is a primitive solvable subgroup, then
an element $g \in G^{\#}$ has at most $2$ fixed points.
\ep
\begin{proof}[Proof of ~\ref{fix1} and \ref{fix2}]
Let $A$ be the unique minimal normal subgroup of $G$.
If $x \in \Omega^g$, then $\Omega = Ax$ and $\Omega^g = C_A(g)x$, where
$C_A (g)$ denotes the centralizer of $g$ in $A$.
Since $\# A = p^k$ and
$C_A (g)$ is a proper subgroup, the conclusion of Proposition~\ref{fix1}
holds.

Now assume $p=2$ and $d-1$ is a prime number. Pick $g \in G^{\#}$
and $x \in \Omega^g$, so that $G_x$ is non trivial.
Let $B \subset G_x$ be a minimal normal subgroup.
Since $G_x$ is solvable, $B$ is an elementary abelian $q$ group for some
prime $q$. We have $\Omega^B = C_A (B) x$.

Note that $G_x$ acts by conjugation on $A$, and, since $G_x$ is a maximal
subgroup, $A$ has no  proper $G_x$-invariant subgroup. Therefore either
$C_A (B)=\{1\}$ or $C_A (B) = A$. If we had $A=C_A (B)$, then
$\Omega^B=\Omega$, that is, $B$ would act
trivially on $\Omega$ and so $B = \{ 1 \}$, which is absurd.
Thus $C_A (B) = \{ 1 \}$ and  $\Omega^B = \{x\}$: $x$
is the only fixed point of the action of $B$ on $\Omega$.

Now let $\Omega_1 = \Omega \smallsetminus \{x\}$.
Since $\# B = q^n$ with $q$ prime and no point of $\Omega_1$ is
fixed by $B$, the cardinality of $\Omega_1$ is a power of $q$.
But by assumption $\# \, \Omega_1= d-1$ is a prime number, hence
$\# \, \Omega_1= q$. It follows immediately that
$B$ and $G_x$ are transitive - hence primitive because $q$ is prime -
subgroups of $\mbox{Aut} (\Omega_1)$. Now we can conclude applying
the first statement to $G_x \subset \mbox{Aut} (\Omega_1)$.
\end{proof}

\section{Zariski's argument} \label{three}
Given a branch point $y$ of a ramified
covering $X \ra Y$, we denote
by $b(y)$ the multiplicity of $y$ in the branch divisor.
Zariski uses~\ref{fix1} to prove the lower bound (\ref{1})
in the following proposition:
\bp \label{zar}
Let $\pi: X \ra Y$ is  a degree $d$ primitive and solvable covering of curves.
Then there exists
a prime $p$ such that $d=p^k$, and for every branch point $y$ we have
\begin{equation} \label{1}
b(y) \geq  \frac{p^k-p^{k-1}}{2} \; .
\end{equation}
If $p=2$ and $d-1$ is a Mersenne prime, then
\begin{equation} \label{mersenne}
b(y) \geq 2^{k-1}-1 \; .
\end{equation}
\ep
\begin{proof}
Let $\tilde{X} \ra Y$ be the Galois closure of the covering $\pi:X \ra Y$.
The Galois group of $\pi$ is by definition the the group
$G= \mbox{Aut} (\tilde{X}/Y) \cong \mbox{Gal} (K(\tilde{X}/K(Y))$
with its action on the generic fiber $\Omega$ of $\pi$.

Since $G \subset \mbox{Aut} (\Omega)$ is solvable and primitive,
by Propositions~\ref{prim} and~\ref{fix1} there is a prime $p$ such that
$d = p^k $ and for every $g \in G^{\#}$ the number
of fixed points of $g$ is at most $p^{k-1}$. If $d=2^k$ and $d-1$ is a Mersenne prime, by
Proposition~\ref{fix2} the number of fixed points of $g$ is at most $2$.

Suppose now $y$ is a branch point of $\pi$, with
$$
\pi^{*} (y) =  \sum_{1}^{t} m_i z_i + z_{t+1} + \cdots + z_{t+n}
$$
where the $z_i's$ are all distinct and $m_i \geq 2$.
Then a small loop around $y$ yields an element $g \in G$
with $n$ fixed points. Since $ n \leq p^{k-1}$ (resp. $n \leq 2$ if $d=2^k$
and $d-1$ is  prime),
we have
$$m=_{def} \sum_{1}^{t} m_i \geq p^k-p^{k-1} \;\;\;\mbox{(resp. $m \geq 2^k-2$)}.$$

Since $m= \sum_{1}^{t} m_i \geq 2t$, we have
$$
b(y) = \sum_{1}^t (m_i-1) = m - t \geq \frac{m}{2} \geq
\frac{p^k-p^{k-1}}{2} \;\;\;\mbox{(resp. $b(y) \geq 2^{k-1}-1$)}.
$$
\end{proof}

Now note that, if
$b(y) \geq l$ for every branch point of a covering $\pi: X \ra Y$,
then by Riemann-Hurwitz:
\begin{equation}\label{due}
lr \leq \deg B_{\pi} = 2 (g(X)-2) - d (2 g(Y) -2)
\end{equation}
Therefore the family of curves $X$ of genus $g$ admitting a $d:1$ covering of $\PP^1$
(resp. of some elliptic curve) 
having $r$ distinct branch points, each with multiplicity at least $l$, has dimension
at most:
\begin{equation} \label{tre}
\frac{1}{l} ( 2 g - 2 + 2 d) - 3 \;
\mbox{  (resp. } \frac{1}{l} ( 2 g - 2) \,).
\end{equation}

Combining this remark with Proposition~\ref{zar} - we need
inequality~(\ref{mersenne}) to deal with the case $d=8$ -
we obtain
\bc \label{cor}
Let $\mathcal{Z}$ be an irreducible family of curves
of genus $7$. Suppose that the general curve in $\mathcal{Z}$ is a
primitive and solvable covering of degree $d \geq 5$ of a curve of
genus $0$ or $1$. Then $\dim \mathcal{Z} \leq 8$.
\ec

\begin{example}
We give another example where inequality~(\ref{mersenne}) is useful.
Consider the $7$-dimensional irreducible family
of all hyperelliptic curves of genus $4$. Let $X$ be a general curve in the family,
and let $\pi: X \ra \PP^1$ be the hyperelliptic projection. We claim that all
solvable maps $f: X \ra \PP^1$ factor through $\pi$. Indeed, a dimension count
based on Riemann-Hurwitz shows the general curve $X$ has no map onto a curve
$Y$ of genus $>0$.
On the other hand, (\ref{three}) and Proposition~\ref{zar} show $X$
has no primitive solvable map to $\PP^1$ of degree $d \geq 5$. Finally, by Riemann-Roch,
a map $f: X \ra \PP^1$ of degree $\leq 4$ is defined by a special linear series, hence
must factor through the hyperelliptic projection $\pi: X \ra \PP^1$.
\end{example}

\section{Debarre-Fahlaoui's construction} \label{four}
We briefly review the construction of~\cite{df}. We will use the following notation:
\begin{itemize}
  \item $E$ is an elliptic curve with origin $o$; if $x$ is a point of $E$,
  we denote by $(x)$ the corresponding divisor.
  \item $S = E^{(2)} =$ second symmetric product of $E$; we will identify
  its closed points with degree two effective divisors $(x)+(y)$ on $E$.
  \item The Abel sum map $s: (x) + (y) \mapsto x+y$ makes $S$ into a ruled
  surface over $E$. We denote by $F$ the fiber over the origin and by
  $C$ the section $\{(o) + (x) : \,\, x \in E \}$.
  \item
  We have $S \cong \PP(\cae)$ where $\cae$ is the rank two bundle over $E$ that fits
  into a nontrivial exact sequence:
$$
\exact{\coo_E}{\cae}{\coo_E (o)}.
$$
Under this identification,
the section $C \hookrightarrow S$ corresponds to the surjection
$\cae \ra \coo_E (o)$. The Euler sequence on $\PP(\cae)$ shows
$\omega_S \cong \coo_S (-2C + F)$, that is, $-K= 2C-F$ is the anticanonical
divisor. Note that $-2K = \Delta$ where
$$\Delta = \{(x)+(x) : \;\; x \in E \}$$
is the diagonal, which is the branch curve of the double cover
$E \times E \ra E^{(2)}=S$.
\item
The divisor $H=3C-K=5C-F$ is very ample, and we denote by $Z$ a smooth irreducible
curve in the linear system $H$. We have $g(Z)=7$ and $\dim \: |H| = 8$.
\end{itemize}

\bp \label{moduli}
The image of the moduli map $|H| \ra \mathcal{M}_7$ is $8$ dimensional.
Therefore varying the elliptic curve $E$, the above construction
yields a $9$ dimensional family $\mathcal{Z}$ of genus $7$ curves $Z$.
\ep
\begin{proof}
Fix $Z$ as above. It is enough to show that the differential of the moduli
map is injective at $Z$. This is a standard technique, but we include some details
for the convenience of the reader.

The relative sheaf of differentials $\Omega_{S/E}$ is isomorphic to
$\omega_s$. Since the anticanonical divisor
has no sections, the exact sequence of tangent sheaves:
$$
\exact{\omega_{S}^{-1}}{\mathcal{T}_S}{s^*(\mathcal{T}_E)}
$$
splits, hence $\mathcal{T}_S \cong \omega_S^{-1} \oplus \coo_S$ and
in particular $h^0 (\mathcal{T}_S) = 1$ and
$h^0(\mathcal{T}_S (-Z))= h^1(\mathcal{T}_S (-Z))= 0$.
Therefore
$h^0(\mathcal{T}_S |_{Z} )=1$.

The long exact cohomology sequence associated to
$$
\exact{\mathcal{T}_Z}{\mathcal{T}_S|_Z}{\mathcal{O}_Z (Z)}
$$
implies that the kernel of
$$
H^0 (\coo_Z (Z))\stackrel{g}{\ra} H^1 ( \mathcal{T}_{Z} )
$$
is one dimensional. Now $g$ is the differential
at $Z$ of the moduli map from the Hilbert scheme of curves algebraically
equivalent to $Z$. The kernel $\mathbb{C} = H^0 (\mathcal{T}_S|_Z)$
is the tangent line to the deformations of $Z$ given by the translations
$\tau_a:(x)+(x) \mapsto (x+a)+(x+a)$
therefore the composite map
$$\mathbb{C} = H^0 (\mathcal{T}_S|_Z) \ra H^0 (\coo_Z (Z)) \ra H^1 ( \coo_S )$$
is an isomorphism.
Thus the differential of $|H| \ra \mathcal{M}_7$ is injective at $Z$.
\end{proof}

\section{Conclusion} \label{five}
\bp \label{hyp}
Let $Y$ be a curve algebraically uniformized by radicals. If
there is a nonconstant morphism from $Y$ to the $4th$ symmetric product
of a curve $X$, then $X$ is also algebraically
uniformized by radicals.
\ep
\begin{proof}
We may assume that $Y$ is rationally uniformized by radicals.
Let $C \subset X \times Y$ be the curve of the correspondence associated to the
map $Y \ra X^{(4)}$.
At least one irreducible component $C_0$ of $C$ projects surjectively both onto $X$ and
onto $Y$. By construction, the normalization $\widetilde{C_0}$ of $C_0$ has a map of degree
at most $4$ onto $Y$, hence $\widetilde{C_0}$ is rationally uniformized by radicals.
Since $\widetilde{C_0}$ covers $X$, we are done.
\end{proof}

\bt
Let $Z$ be a general curve in the family $\mathcal{Z}$ of Proposition~\ref{moduli}.
Then $Z$ is algebraically but not rationally uniformized by radicals.
\et
\begin{proof}
We first show $Z$ is algebraically uniformized by radicals.
For every $p \in E$ let $C_p$ denote the curve $\{(p) + (x) : \,\, x \in E \}$ on $S=E^{(2)}$. Then
$Z$ and $C_p$ intersect in a divisor of degree $4$. This gives rise to a nonconstant
morphism from $E$ to the 4th symmetric product of $Z$, and by Proposition~\ref{hyp}
$Z$ is algebraically uniformized by radicals.

Now we prove $Z$ is not rationally uniformized by radicals. For this it is enough
to show there is no nontrivial solvable primitive covering $f: Z \ra W$, where $W$ is
an arbitrary curve.

As in~\cite{df}, we note that $Z$ does not cover any curve $W$ of genus $> 1$ because
the quotient of the Jacobian $JC$ of $C$ by $s^*( JE)$ is simple; for the same reason
$Z$ does not have a morphism of degree $d \leq 4$ onto an elliptic
curve~\cite[Proposition 5.14]{df}.

By ~\cite{d}, the curve $Z$ does not have a morphism of degree $d \leq 4$
onto a rational curve either.

Thus we only have to exclude primitive and solvable coverings $f: Z \ra W$
of degree $d \geq 5$ where $W$ is a rational or an elliptic curve. This is taken
care by Corollary~\ref{cor}, because $Z$ varies in a family of dimension $9$ by
Proposition~\ref{moduli}.
\end{proof}

\end{document}